\documentclass[12pt]{article}
\usepackage{amssymb,latexsym}
\usepackage{fullpage}
\usepackage{array}
\renewcommand{\thefootnote}{}

\begin{document}
\title{Cohomology of toric bundles}
\author{P.Sankaran and V.Uma    \\  
Institute of Mathematical Sciences\\
CIT Campus, Chennai 600 113, INDIA \\
E-mail: {\tt sankaran@imsc.res.in}\\ 
 uma@imsc.res.in\\[4mm]} 
\date{}
\maketitle

\footnote{2000 A.M.S. Subject Classification:- 14M25, 14F15\\
Key words and phrases: Toric varieties, toric bundles, singular cohomology, 
Chow ring, and K-theory.\\
This is the corrected version of the paper published in Comm. Math. Helv. 
{\bf 78}, (2003), 540--554. 
The {\it Errata} appeared in Comm. Math. Helv. {\bf 79}, (2004), 840--841. }

\thispagestyle{empty}

\def\theequation {\arabic{section}.\arabic{equation}}
\renewcommand{\thefootnote}{}

\newcommand{\codim}{\mbox{{\rm codim}$\,$}}
\newcommand{\stab}{\mbox{{\rm stab}$\,$}}
\newcommand{\lr}{\mbox{$\longrightarrow$}}

\newcommand{\ch}{{\cal H}}
\newcommand{\cf}{{\cal F}}
\newcommand{\cd}{{\cal D}}

\newcommand{\blr}{\Big \longrightarrow}
\newcommand{\da}{\Big \downarrow}
\newcommand{\ua}{\Big \uparrow}
\newcommand{\hra}{\mbox{\LARGE{$\hookrightarrow$}}}
\newcommand{\rt}{\mbox{\Large{$\rightarrowtail$}}}
\newcommand{\dua}{\begin{array}[t]{c}
\Big\uparrow \\ [-4mm]
\scriptscriptstyle \wedge \end{array}}

\newcommand{\be}{\begin{equation}}
\newcommand{\ee}{\end{equation}}

\newtheorem{guess}{Theorem}[section]
\newcommand{\bth}{\begin{guess}$\!\!\!${\bf .}~}
\newcommand{\eeth}{\end{guess}}
\renewcommand{\bar}{\overline}
\newtheorem{propo}[guess]{Proposition}
\newcommand{\bpropo}{\begin{propo}$\!\!\!${\bf .}~}
\newcommand{\epropo}{\end{propo}}

\newtheorem{lema}[guess]{Lemma}
\newcommand{\blem}{\begin{lema}$\!\!\!${\bf .}~}
\newcommand{\elem}{\end{lema}}

\newtheorem{defe}[guess]{Definition}
\newcommand{\bdefe}{\begin{defe}$\!\!\!${\bf .}~}
\newcommand{\edefe}{\end{defe}}

\newtheorem{coro}[guess]{Corollary}
\newcommand{\bcor}{\begin{coro}$\!\!\!${\bf .}~}
\newcommand{\ecor}{\end{coro}}

\newtheorem{rema}[guess]{Remark}
\newcommand{\brem}{\begin{rema}$\!\!\!${\bf .}~\rm}
\newcommand{\erem}{\end{rema}}

\newtheorem{exam}[guess]{Example}
\newcommand{\beg}{\begin{exam}$\!\!\!${\bf .}~\rm}
\newcommand{\eeg}{\end{exam}}

\newcommand{\ctext}[1]{\makebox(0,0){#1}}
\setlength{\unitlength}{0.1mm}
\newcommand{\cl}{{\cal L}}
\newcommand{\cp}{{\cal P}}
\newcommand{\ci}{{\cal I}}
\newcommand{\bz}{\mathbb{Z}}
\newcommand{\cs}{{\cal s}}
\newcommand{\cv}{{\cal V}}
\newcommand{\ce}{{\cal E}}
\newcommand{\ck}{{\cal K}}
\newcommand{\cR}{{\cal R}}
\newcommand{\bq}{\mathbb{Q}}
\newcommand{\bt}{\mathbb{T}}
\newcommand{\bh}{\mathbb{H}}
\newcommand{\br}{\mathbb{R}}
\newcommand{\wt}{\widetilde}
\newcommand{\im}{{\rm Im}\,}
\newcommand{\bc}{\mathbb{C}}
\newcommand{\bp}{\mathbb{P}}
\newcommand{\spin}{{\rm Spin}\,}
\newcommand{\ds}{\displaystyle}
\newcommand{\tor}{{\rm Tor}\,}
\newcommand{\bs}{\mathbb{S}}
\def\ns{\mathop{\lr}}
\def\nssup{\mathop{\lr\,sup}}
\def\nsinf{\mathop{\lr\,inf}}
\renewcommand{\phi}{\varphi}
\newcommand{\co}{{\cal O}}

\noindent
{\bf Abstract:} 
Let $p:E\lr B$ be a principal bundle with 
fibre and structure 
group the torus $T\cong (\bc^*)^n$ over a topological space $B$. 
Let $X$ be a nonsingular projective $T$-toric variety. 
One has the $X$-bundle $\pi:E(X)\lr B$ where $E(X)=E\times_T X$, $\pi([e,x])
=p(e)$.  This  is a Zariski locally trivial fibre bundle 
in case $p:E\lr B$ is algebraic. The purpose of this 
note is to describe (i) the singular cohomology ring of $E(X)$ as 
an $H^*(B;\bz)$-algebra, 
(ii) the topological K-ring of 
$K^*(E(X))$ as a $K^*(B)$-algebra when $B$ is compact Hausdorff. 
When $p:E\lr B$ is algebraic over an irreducible, nonsingular, 
noetherian scheme over $\bc$, we describe  
(iii) the Chow ring of $A^*(E(X))$ as an $A^*(B)$-algebra, and  
(iv) the Grothendieck ring $\ck^0(E(X))$ of algebraic vector 
bundles on $E(X)$ as a $\ck^0(B)$-algebra. 

\section{Introduction}
Let $T\cong (\bc^*)^n$ denote the complex algebraic torus.  Let $M=Hom_{alg}
(T,\bc^*)\cong \bz^n$ and $N=Hom_{alg}(\bc^*,T)\cong \bz^n$ denote the group 
of characters and the group of $1$-parameter subgroups of $T$ respectively.  
Note that $M=N^\vee:=Hom(N;\bz)$ under the natural pairing $
\langle~,~\rangle :M\times N\lr \bz,$ given by $\chi^u\circ \lambda_v(z)=z^
{\langle u,v\rangle} $ for all $z\in \bc^*$. (Here $\chi^u\in \bc(T)=\bc[M]$
denotes the character corresponding 
to $u\in M$ and $\lambda_v$ the $1$-parameter subgroup corresponding to 
$v\in N$.)  

Let $\Delta$ be a fan in $N$ such  that the $T$-toric variety 
$X:=X(\Delta)$ is complete and non-singular.  
Let $p:E\lr B$ be a principal bundle with structure group and 
fibre the torus $T$ over an arbitrary topological space $B$. 
When the  bundle $p:E\lr B$ is algebraic, 
it is well-known that the bundle $E\lr B$ is Zariski locally trivial. 

Consider the fibre  bundle $\pi:E(X)\lr B$ with fibre the 
toric variety $X$, 
where $E(X)$ is the fibre product $E\times_T X$, and the projection map 
is defined as $\pi([e,x])=p(e)$. Note that the bundle $\pi:E(X)\lr B$ 
is Zariski locally trivial when $p:E\lr B$ is algebraic.  
In this paper, we describe the integral singular cohomology ring 
$H^*(E(X);\bz)$, and the $K$-ring $K(E(X))$ when $B$ is a compact topological 
space. Also, when $p:E\lr B$ is algebraic and $B$ an irreducible  nonsingular 
noetherian scheme over $\bc$, we describe the Chow ring $A^*(E(X))$, and the 
Grothendieck ring $\ck^0(E(X))$ of algebraic vector bundles of 
the complex variety $E(X)$.  

Suppose that $q:V\lr X$ is a $T$-equivariant vector bundle over $X$, then 
we obtain a vector bundle $E(V)$ over $E(X)$ with total space $E\times_T V$ 
where the bundle projection is the map $[e,v]\mapsto [e,q(v)]$.  
In case $V$ is a $T$-equivariant line bundle  
associated to a character $\chi^u:T\lr \bc^*,$ the bundle 
$E(V)$ is isomorphic to the pull-back 
bundle $\pi^*(\xi_u)$ where $\xi_u\lr B$ is the line bundle got from 
$E\lr B$ by `extending' the structure group via $\chi^u$. 
After fixing an isomorphism, 
$T\cong (\bc^*)^n$,  $u\in M$ corresponds to an element $(a_1,\cdots, a_n)\in 
\bz^n$.  The bundle $E$ is then the principal bundle associated to 
the Whitney sum of line bundles $\xi_i, 1\leq i\leq n,$ and  
$E(V)$ can then be identified with the tensor product 
$\xi_1^{a_1} \otimes \cdots \otimes \xi_n^{a_n}$.  
(Here it is understood that, when $a<0$, $\xi^a=(\xi^\vee)^{-a}$.)  

When $B$ is a  non-singular variety 
any line bundle $\xi$ over $B$ is 
isomorphic to $\co(Y)$ for some divisor $Y$ in $B$. The divisor  
class $[Y]$ is the first Chern class $c_1(\xi)\in A^1(B)$  of $\xi$.   

We use the notations of \cite{f} throughout the paper.

For $k \geq 1$, $\Delta(k)$ will denote the set of $k$ dimensional cones in 
$\Delta$. We let $d=\#\Delta(1)$, and write $v_1,\cdots, v_d$ for the 
primitive elements of $N$ along the edges in $\Delta$. Let  
$\rho_j\in \Delta(1)$ be the edge $\br_{\geq 0}v_j$. Recall that 
our hypothesis that $X$ is smooth is equivalent to the statement that the 
set of the primitive vectors along the edges of any cone in $\Delta$ 
is part of a $\bz$-basis for $N$. 

For a cone $\sigma\in \Delta$, $U_\sigma$ denotes the affine toric variety 
defined by $\sigma$ and $V(\sigma)$ denotes the closure in $X$ of 
the variety whose local equation in $U_\sigma$ is $\chi^u=0$ for all 
$u\notin \sigma^\perp$, $u\in \sigma^\vee$. The $V(\sigma),\sigma\in\Delta,$ 
are the orbit closures for the action of $T$ on $X$. 
  
For $1\leq j\leq d$ let $L_j$ denote the $T$-equivariant line 
bundle over $X$ which 
corresponds to the piecewise linear function $\psi_j$ defined by $\psi_j(v_i)
=-\delta_{i,j}$. The line bundle $L_j$ admits a global $T$-equivariant 
section $s_j$ whose zero locus is the variety $V(\rho_j)$. 

Let $\sigma_1,\cdots,\sigma_m$ be an ordering of the cones in $\Delta(n)$.
Let $\tau_i\in \Delta$ be the intersection with $\sigma_i$ of those cones 
$\sigma_j$, $j>i,$ such that $\dim(\sigma_i\cap\sigma_j)=n-1.$  
Thus $\tau_1=0$, and $\tau_m=\sigma_m$.
Consider the condition:
$$\tau_i<\sigma_j \Longrightarrow i\leq j. \eqno {(*)}.$$  
Set $\tau'_i<\sigma_i$ to be the cone such that $\tau_i\cap\tau'_i=0, 
\dim(\tau_i)+\dim(\tau'_i)=n$, $1\leq i\leq m$. Also consider the condition 
$$\tau_i'<\sigma_j\Longrightarrow j\leq i. \eqno {(*')}$$ 
Note that $\tau'_i$ is the intersection with $\sigma_i$ of those cones 
$\sigma_j$ with $j<i$ and $\dim(\sigma_i)\cap\dim(\sigma_j)=n-1$ and so 
condition $(*')$ is the same as $(*)$ for the reverse ordering 
on $\Delta$. 
It is well known that when $X$ is (nonsingular) projective, then 
there exists an ordering of the cones in $\Delta(n)$ such that both 
conditions  $(*)$ and $(*')$  hold.  We shall assume that 
there exists an ordering of $\Delta(n)$ such that property $(*)$ holds. 
(See \cite{f}, \S 5.2.)  

By relabelling the $v_j$'s if necessary, we assume that 
$v_1,\cdots, v_n \in N$ are primitive vectors along the edges of 
$\sigma_m$ and let $u_1,\cdots, u_n$ be the dual basis 
of $M$. 

\bdefe
Let $S$ be a ring with $1$.  
Let $r_i, 1\leq i\leq n,$ be in the centre of $S$. Consider 
the polynomial algebra $S[x_1,\cdots, x_d]$.  We denote by $I_S$ the two-sided 
ideal generated by the following two types of elements 
$$x_{j_1}\cdots x_{j_k}, ~1\leq j_p\leq d, \eqno{(i)} $$ 
where $v_{j_1},\cdots,v_{j_k}$ do not span a cone of $\Delta$, 
and,  
$$y_i:= \sum_{1\leq j\leq d} \langle u_i,v_j\rangle x_j-r_i, 
~1\leq i\leq n. \eqno{(ii)}$$ 
\noindent
Denote by $\ci_S$ the two-sided ideal 
generated by elements of type $(i)$  above and the elements  
$$z(u):=\prod_{j, \langle u,v_j\rangle>0} 
(1-x_j)^{\langle u,v_{j}\rangle}-
r(u)\prod_{j, \langle u,v_{j}\rangle<0}(1-x_j)^{-\langle u,v_{j}\rangle}, 
~u\in M, \eqno{(ii')}$$  
where $r(u)=\prod_{1\leq j\leq d}r_i^{\langle u,v_j\rangle}.$
Define $R=R(S,\Delta)=S[x_1,\cdots, x_d]/I_S$ and $\cR=\cR(S,\Delta) 
=S[x_1,\cdots,x_d]/\ci_S$.
\edefe

Note that the $S$-algebras $R$ and $\cR$ depend not only on the fan $\Delta$, 
but also on the  
the isomorphism $N\cong \bz^n$ resulting from the choice of 
$\sigma_m\in \Delta$ and the elements $r_i\in S$.  The only non-commutative 
ring $S$ we need to consider is the integral cohomology ring of $B$.  

Note that for any cohomology theory $\ch$, 
$\ch^*(E(X))$ is an $\ch^*(B)$-algebra via the induced map $\pi^*:\ch^*(B)
\lr \ch^*(E(X))$.  The following is our main theorem:

\noindent
\bth \label{main} Let $\pi:E\lr B$ be a principal $T$-bundle 
over an arbitray topological space $B$. 
Assume that $X$ is a smooth complete $T$-toric variety 
and that $\Delta(n)$ has been ordered so that $(*)$ holds.  
With above notations,\\ 
(i)  The singular cohomology ring of $E(X)$ is isomorphic as an 
$H^*(B;\bz)$-algebra to $R(H^*(B;\bz),\Delta)$, with  
$r_i=c_1(\xi_i^\vee)\in H^2(B;\bz)$. \\
(ii) When $B$ is compact Hausdorff, 
the $K$-ring $K^*(E(X))$ of complex vector bundles 
over $E(X)$ is isomorphic as a $K^*(B)$-algebra to $\cR(K^*(B);\Delta)$ where 
$r_i=[\xi_i^\vee]\in K(B)$, $1\leq i\leq n$.\\ 
Suppose  $p:E\lr B$ is algebraic where $B$
irreducible, non-singular and noetherian over $\bc$. Furthermore, assume that  
$(*')$ also holds. Then:\\
(iii) The Chow ring $A^*(E(X))$ of $E(X)$ is isomorphic as an $A^*(B)$-algebra 
to $R(A^*(B),\Delta)$ where 
$r_i=c_1(\xi_i^\vee)\in A^1(B)$, $1\leq i\leq n$.\\
(iv) The ring $\ck(E(X))$ is isomorphic as a $\ck(B)$-algebra  
to $\cR(\ck(B),\Delta)$ where $r_i=[\xi_i^\vee]\in \ck(B)$.  
\eeth

We shall now briefly explain the method of proof. For the first three 
parts, we shall use a Leray-Hirsch type theorem to obtain the 
structure of $\ch^*(E(X))$ as a {\it module} over $\ch^*(B)$. 
Then we shall construct a $\ch^*(B)$-algebra  homomorphism from the 
expected $\ch^*(B)$-algebra to $\ch^*(E(X))$ and verify that this 
algebra homomorphism is an isomorphism of $\ch^*(B)$ -{\it modules}. 
The ``Leray-Hirsch'' theorem  we need in the context of Chow rings 
is due to D.Edidin and W.Graham \cite{eg}. However we give a proof which 
is more suited to our specific situation.
(See also \cite{es}.)  
The ``Leray-Hirsch'' in the context of $K$-theory of complex vector bundles 
that we need is theorem 1.3, Ch. IV, \cite{k}. 
(Compare Theorem 2.7.8, \cite{at}.)  For part (iv) we 
use a result of Grothendieck \cite{sga} to prove the analogue of  
Leray-Hirsch theorem.  

We do not know if parts (iii) and (iv) of  the main theorem remain valid 
without the hypothesis that  $(*')$ hold. Neither do we know of an 
example where $\Delta(n)$ admits an ordering  satisfying $(*)$ but no 
ordering  
that satisfies both $(*)$ and $(*')$. However, there are complete nonsingular  
varieties $X(\Delta)$ which are {\it not} projective such that $\Delta(n)$ 
admits an ordering satisfying both $(*)$ and $(*')$. The  example of  
a complete non projective toric variety given in p. 84 \cite{oda} is 
seen to be one such.  

Examples of algebraic bundles $E(X)\lr B$ we consider include as special 
cases the toric fibre bundles considered on p.41, \cite{f}. 

We were motivated by the work of Al Amrani \cite{ala} who has computed 
the singular cohomology ring of a weighted projective space bundle. 
Another motivation for us was  the work of H.Pittie and A.Ram \cite{pr} who 
established the Pieri-Chevalley formula in $K$-theory 
in the context of an algberaic $G/B$ 
bundle associated to a principal $B$ bundle where $G$ is a complex 
simple algebraic group and $B$ a Borel subgroup.

\section{The rings $R$ and $\cR$}

In this section we prove certain facts about the rings $R$ and $\cR$ 
which will be needed in the proof of the main theorem.

We keep the notations of \S1.  We assume that $\Delta(n)$ has been 
so ordered that property $(*)$ holds.
Recall that $v_1,\cdots,v_d$ are the primitive 
vectors along the edges of $\Delta$, that $v_1,\cdots, v_n$ 
are in $\sigma_m$, and that $u_1,\cdots, u_n$ is the dual basis of 
$M$.  

For any cone $\gamma\in \Delta$, denote by $x(\gamma)$ the monomial 
$x_{j_1}\cdots x_{j_r}\in S[x_1,\cdots,x_d]$  
where $v_{j_1}\cdots ,v_{j_r}$ are 
the primitive vectors along the edges of $\gamma$.  

Recall from \S1 the definition of the $S$-algebras $R$ and $\cR$. 
We shall denote by the same symbol 
$x(\gamma)$, in $R$ and $\cR$, the image of the monomial 
$x(\gamma)\in S[x_1, \cdots, x_d]$ under the canonical quotient map.

\noindent
\blem \label{rmod} 
\noindent
(i) For $u=\sum_{1\leq i\leq n}a_iu_i\in M$, the following equality  
holds in $R$. 
$$\sum_{1\leq j\leq d}\langle u, v_j\rangle x_j=r_u$$
where $r_u=\sum_{1\leq i\leq n} a_ir_i.$  \\
\noindent 
(ii) If $\gamma\in\Delta(r)$ is spanned 
by $v_{j_1}, v_{j_2},\cdots, v_{j_r}$, 
then 
$$x_{j_1}x(\gamma)=-\sum_k \langle u,v_k\rangle x(\gamma_k) +r_ux(\gamma)$$
for some $u\in M$, 
where the sum on the right is over those cones $\gamma_k $ in $\Delta(r+1)$ 
which are spanned by primitive vectors 
$v_k,v_{j_1}, v_{j_2},\cdots, v_{j_r}$.\\
\noindent 
(iii) If $\alpha<\gamma\leq \beta$ are cones in $\Delta$ then there 
exist cones 
$\gamma_1,\cdots,\gamma_s\in \Delta$ with $\alpha<\gamma_k$ such that the 
$\gamma_k$ are not contained in $\beta$, and   
$$x(\gamma)=\sum_{k}c_kx(\gamma_k)+c x(\alpha)$$ for some $c, c_k\in S$.\\ 
\noindent
(iv) The monomials $x(\tau_i)$, $1\leq i\leq m,$ span  $R$ as an $S$-module.  
\elem

\noindent
{\bf Proof:} Part (i) follows easily from linearity of the pairing 
$\langle ~,~\rangle$ with respect to the first argument. \\

\noindent
(ii) Suppose $\gamma\leq \sigma$ where $\sigma$ is $n$-dimensional. 
Let $v_{j_1}, \cdots,v_{j_n}$ be the primitive vectors which span 
$\sigma$ such that the first $r$ elements span $\gamma$, with $j_1=j$.  
Let $u\in M$ be the dual basis element such that $\langle u,v_{j_k}
\rangle =\delta_{j,j_k}$.
One has the relation:
$$x_{j_1}+\sum_{k\neq j_1}\langle u,v_k\rangle x_k-r_u=0.$$ 
Multiplying both sides by $x(\gamma)$, and using the type $1$  
relations, we get   
$$x_{j_1}x(\gamma)=-\sum_k \langle u,v_k\rangle x(\gamma_k) +r_ux(\gamma)$$
where the sum on the right is over those cones $\gamma_k $ in $\Delta(r+1)$ 
which are spanned by primitive the vectors $v_k,v_{j_1}v_{j_2},\cdots, 
v_{j_r}$ where  $k\neq{j_1},\cdots,{j_n}$.
This proves (ii).

\noindent 
(iii) Suppose $v_{j_1},\cdots, v_{j_l}$ spans $\beta\in \Delta(l)$ such 
that the 
first $r$ of these span $\alpha$ and the first $p$ of these span $\gamma$, 
$p>r$. Without loss of generality we may assume that 
$\beta$ is an $n$ dimensional cone so that 
$l=n$, and $v_{j_1},\cdots, v_{j_n}$ is a basis for $N$.  Now let $u\in M$ be 
the dual basis element so that $\langle u,v_{j_q}\rangle=\delta_{p,q}$. 
Then we have 
$$x_{j_p}+\sum_{k\neq j_p} \langle u,v_k\rangle x_k -r_u=0. \eqno(1)$$ 
Multiplying by $x_{j_1}\cdots x_{j_{p-1}} $ and observing that 
the coefficient of $x_k$ in the sum is zero for $k\in\{j_1,\cdots, j_n\}$ 
and $k\neq j_p $, we get 
$x(\gamma)+\sum \langle u,v_k\rangle x(\gamma_k)-r_ux(\gamma')=0$ where 
$\gamma'$ is the cone spanned by $v_{j_1},\cdots,v_{j_{p-1}}$ and 
the sum is over those cones $\gamma_k\in \Delta(p)$ which are spanned by 
$v_{j_1},v_{j_2},\cdots,v_{j_{p-1}},v_k,~ k\neq {j_1,\cdots,j_n}$. Note 
that each of these $\gamma_k$ contains $\alpha$ but is not contained in 
$\beta$. If $\gamma'=\alpha$, we are done. Otherwise, 
by an induction on the dimension of $\gamma$ 
the statement is true for $\gamma'$.  
Substituting this expression for $x(\gamma')$ in (1), we see that (iii) holds. 

\noindent 
(iv) We first prove that the $x(\tau_i)$ span $R$.  
In view of (ii), it suffices to prove that for any $\gamma$, 
$x(\gamma)$ is in the $S$-submodule spanned by the $x(\tau_i)$.  
Property $(*)$ implies that given any $\gamma\in \Delta$, there exists 
a unique $i$ such that $\tau_i\leq \gamma\leq \sigma_i$; indeed it is 
the smallest $i$ for which $\gamma\leq \sigma_i$. (See \cite{f}, \S5.2.)  
We prove, by a downward induction 
on this $i$, that $x(\gamma)$ is in the $S$-span of  $x(\tau_j),~j\geq i$.  
If $i=m$, then $\gamma=\sigma_m=\tau_m$ and there is nothing to prove.

Let $\tau_i\leq \gamma\leq \sigma_i$ for some $i<m$. Now, using (iii), we 
can write 
$x(\gamma)$ as an $S$-linear combination of $x(\tau_i)$ and 
$x(\gamma_j)$ where 
$\tau_i<\gamma_j$, and $\gamma_j$ is not contained in $\sigma_i$.  
It follows that each $\gamma_j$ is such that $\tau_{r}\leq \gamma_j\leq 
\sigma_{r}$ for some $r$ (depending on $j$) with $r>i$. By inductive 
hypothesis, each of the $x(\gamma_j)$ is in the $S$-span of $x(\tau_q), 
~q\geq r$.  It follows that  $x(\gamma)$  is in the 
$S$-span of $x(\tau_r), r\geq i$, completing the proof. \hfill $\Box$

Concerning the structure of $\cR$ we have the following.
\noindent 
\blem \label{kring} 
Assume that the elements $r_i\in S, ~1\leq i\leq n$ are invertible.  Then \\
\noindent
(i) For $u=\sum_{1\leq i\leq n} a_iu_i\in M$ the following relation holds 
in $\cR$: 
$$z_u:=\prod_{j, \langle u,v_j\rangle >0}(1-x_j)^{ \langle u,v_j\rangle} 
-r_u\prod_{j, \langle u,v_j\rangle <0}(1-x_j)^{- \langle u,v_j\rangle}=0$$
where $r_u=\prod_{1\leq i\leq n}r_i^{a_i}$.\\
\noindent
(ii)  Let $\alpha<\gamma\leq \beta$ be cones in $\Delta$. Suppose that 
$\gamma$ is spanned by $v_{j_1},\cdots,v_{j_k}$, then 
$$x_{j_1}x(\gamma)=(1-r_u)x(\gamma)+\sum_{p} a_px(\gamma_p)$$ 
where $a_p\in S$, and $\gamma_p\in \Delta$ are such that $\alpha<\gamma_p, 
~\gamma_p$ are not contained in $\beta$ and 
$\dim(\gamma_p)>\dim(\gamma)$. \\
\noindent
(iii) Let $\alpha<\gamma\leq \beta$ be cones in $\Delta$.  Then 
$$x(\gamma)=\sum_p b_px(\gamma_p)+bx(\alpha)$$ 
for some $b_p,b\in S$ and suitable cones 
$\gamma_p\in \Delta$ which contain $\alpha$ and are 
not contained in $\beta$.
\noindent
(iv) The monomials $x(\tau_i), 1\leq i\leq m$ span $\cR$ as an $S$-module.
\elem

{\bf Proof:} Proof of (i) is an easy exercise.\\
(ii) Without loss of generality, we may assume that $\beta$ is an 
$n$-dimensional cone. We prove this by descending induction on 
the dimension of $\gamma$. 
Suppose that $v_{j_1}, \cdots,v_{j_n}$ span 
$\beta$, and that $v_{j_1}\notin \alpha$. Let $u\in M$ be the dual 
basis element  such that $\langle u,v_{j_r}\rangle =\delta_{1,r}$.
The relation $z_u=0$ can be rewritten as 
$$(1-x_{j_1})\prod_{j, \langle u,v_j\rangle>0}(1-x_j)
^{\langle u,v_j\rangle} 
=r_u\prod_{j,\langle u,v_j\rangle<0}(1-x_j)^{-\langle u,v_j\rangle}$$
Note that none of the $x_{j_r}, 2\leq r\leq n$ occur in the 
above relation.  Multiplying both sides by $x(\gamma)$, 
$$(x(\gamma)-x_{j_1}x(\gamma))\prod_p (1-x_p)^{\langle u,v_p\rangle} 
=r_ux(\gamma)\prod_q(1-x_q)^{-\langle u,v_q\rangle} \eqno{(2)} $$ 
where the product is over those $p$, (resp. $q$)  such that $v_p$, 
(resp. $v_q$),  $v_{j_1}\cdots,v_{j_k}$ 
span a cone of $\Delta$, $\langle u,v_p\rangle >0,$ (resp. 
$\langle u, v_q\rangle <0$).
In particular, if $\gamma$ is $n$-dimensional, then the above equation 
reads $x_{j_1}x(\gamma)=(1-r_u)x(\gamma)$, which proves the lemma in this 
case.  Assume that $k<n$ and that the statement holds for all higher 
dimensional cones. Then from equation (2), we see that 
the lemma follows by repeated application of the inductive hypothesis and 
by the observation that if $\gamma'<\gamma''$ and if $\gamma'$ is 
not contained in $\beta$, then neither is $\gamma''$.\\
\noindent
Parts (iii) and (iv) follow from (ii) just as in the proof of lemma \ref{rmod}.
\hfill $\Box$ 

\brem  One can show that if $r_i=0$ for all $i, 1\leq i\leq n$, then 
$x_j^{n+1}=0, ~1\leq j\leq d$ in $R$ and that 
$x(\tau_i), 1\leq i\leq m,$ form a basis for $R$ as a module over $S$. 
Similarly, if $r_i=1$ for all 
$1\leq i\leq n$, then $x_j^{n+1}=0$ for $1\leq j\leq d$ and $x(\tau_i), 
1\leq i\leq m$ form a basis for $\cR$ as an $S$-module.  
\erem

\section{Singular cohomology and Chow ring} 
In this section we shall prove parts (i) and (iii) of the 
Main Theorem \ref{main}.

Let $\Delta$ be a complete nonsingular fan in $N$.
We assume that $\sigma_1,\cdots , \sigma_m$ is an ordering of $\Delta(n)$ 
such that property $(*)$ holds.   (See \S1.)  
This implies that the toric variety $X=X(\Delta)$ has an algebraic 
cell decomposition, namely, there exists closed subvarieties $X=Z_1\supset 
\cdots \supset Z_m$ of $X$ such that 
$Z_i\setminus Z_{i+1}=:Y_i\cong \bc^{k_i}$ 
for some integers $k_i$. In fact, with $\tau_i$ as in \S1,  
the closure of $Y_i$ is just the variety $V(\tau_i)$. See \cite{f}. 
This yields the structure of a (finite) CW complex on $X$ with cells 
only in even dimensions. 

\noindent
{\it Notation:} We shall denote $V(\tau_i)$ by $V_i$.  If $(*')$ also 
holds, then we set $V'_i=V(\tau_i')$.

Assume that $p:E\lr B$ is complex algebraic and 
$B$ irreducible, nonsingular, and  noetherian over $\bc$. 
Now since the varieties $V_i$
are stable under the $T$-action, one has the associated bundles 
$\pi_i:E(V_i)\lr B$  with fibre $V_i$. Note that $E(V_i)$ is a 
smooth closed subvariety of $E(X)$. 
For any closed subvariety $Z$ in an algebraic variety $Y$ we denote by 
$[Z]$ its rational equivalence class in $A_*(Y)$. If $Z$ and $Y$ are smooth,   
we denote by $[Z]$ the cohomology class dual to $Z$ in $H^{2r}(Y;\bz)$ 
as well as the element in the Chow cohomology group $A^r(Y)$ where 
$r$ is the codimension of $Z$ in $Y$.

In case property $(*')$ also holds, then 
$[V_i].[V'_j]=0$ if $j<i$, and 
$[V_i][V'_i]\in H^{2n}(X;\bz)\cong \bz$ is the positive generator 
with respect to the standard orientation coming from the complex structure on 
$X$.  Also, in the Chow ring, $[V_i][V_i']\in A^{n}(X)=A_0(X)\cong 
\bz$ denotes the class of the point $[V(\sigma_i)]$ which generates $A_0(X)$. 

\noindent
\blem\label{leray} Let $X$ be a complete nonsingular $T$-toric variety  
and suppose that property $(*)$ holds for an ordering of $\Delta(n)$.  
Let $\pi:E\lr B$ be a principal $T$-bundle over any topological space.Then:\\ 
\noindent
(i) The bundle $\pi:E(X)\lr B$ admits a cohomology extension of the fibre in 
singular cohomology with integer coefficients. $H^*(E(X);\bz)$ is isomorphic 
to $H^*(B;\bz)\otimes H^*(X;\bz)$ as an $H^*(B;\bz)$-module.\\
\noindent 
Assume $\pi:E\lr B$ is complex algebraic where $B$ an irreducible, nonsingular 
noetherian variety over $\bc$. Suppose that properties  $(*), (*')$  hold.  
Then:\\
\noindent
(ii) The Chow group $A^*(E(X))$ is isomorphic as an $A^*(B)$-module to 
$A^*(B)\otimes A^*(X)$.
\elem
   
\noindent
{\bf  Proof:} We shall fix a base point $b_0\in B$ and identify $X$ with the 
fibre $\pi^{-1}(b_0)\subset E(X)$. 
(i).  Since the $X$ has a CW decomposition with only even 
dimensional 
cells, its integral cohomology is isomorphic to the free abelian group with 
basis labelled by its cells. Indeed the 
dual cohomology classes $[V_i]\in H^{2l_i}(X;\bz)$, $l_i=\dim(\tau_i),$ 
form a $\bz$-basis for $H^*(X;\bz)$. 

Let $\eta\in \Delta(r)$ spanned by the primitive vectors along 
$v_{j_1},\cdots,v_{j_r}$. Denote by $L(\eta)$ the $T$-equivariant 
bundle $L_{j_1}\oplus\cdots \oplus L_{j_r}$, where the $L_j$ are as 
defined in \S1.
The class $[V(\eta)]\in H^{2r}(X;\bz)$ equals the 
the Chern class $c_r(L(\eta))=c_1(L_{j_i})\cdots c_1(L_{j_r}).$  
The bundle 
$\cl(\eta)=E(L(\eta))$ over $E(X)$ restricts to $L(\eta)$ over $X$. By 
the naturality of Chern classes, $c_r(\cl(\eta))\in H^{2r}(E(X);\bz)$ 
restricts to $c_r(L(\eta))=[V(\eta)]\in H^{2r}(X;\bz)$. 
In particular, it follows that $[V(\tau_i)], 
1\leq i\leq m,$ are in the image of the restriction homomorphism 
$H^*(E(X);\bz)\lr H^*(X;\bz)$.  
The lemma follows by Leray-Hirsch theorem (\cite{s}, p.258). 

\noindent
(ii)  Our proof follows that of Lemma 2.8 \cite{es} closely. 
(See also Lemma 6, \cite{eg}.)  
Clearly the classes $[E(V_i)]\in A^*(E(V_i)), 1\leq i\leq m$ restrict 
to elements of a $\bz$ basis (namely $[V_i]\in A^*(X)$). Consider 
the $A^*(B)$-linear map  $\Phi:A^*(B)\otimes A^*(X)\lr A^*(E(X))$, defined as 
$$\Phi(\sum_{1\leq i\leq m}b_i\otimes [V_i])= \sum_{1\leq i\leq m} 
\pi^*(b_i).[E(V_i)].$$ 
To prove (ii) we show that $\Phi$ is an isomorphism. 
Suppose $\Phi(\sum b_i\otimes [V_j])=0 \in A^*(E(X))$.
Assume that $k$ is the smallest integer such that $b_k\neq 0$. Since 
for $j\geq k$, $V_j$ and $V'_k$ are disjoint unless $j=k$ in which case they 
intersect transversally and $V_j\cap V_j'=V(\sigma_j)$ scheme theoretically.  
We see that $[E(V_j)].[E(V'_k)]=0$ if $j>k$, and, $E(V_j),E(V'_j)$ intersect 
transversally and so 
$E(V_j)\cap E(V'_j)=E(V(\sigma_j))$ scheme theoretically (where the 
subvarieties  are given the reduced scheme structure). Therefore, 
$[E(V_j)].[E(V_j')]=[E(V(\sigma_j))]$.  Note that since  $V(\sigma_j)$ is 
a $T$-fixed point, $E(V(\sigma_j))\cong E/T=B$.  Denote by 
$\pi_j$ the restriction of $\pi:E(X)\lr B$ to $E(V(\sigma_j))$.
Also let $\iota_j $ be the inclusion 
$E(V(\sigma_j))\subset E(X)$. Note that $[E(V(\sigma_j))]=
\iota_{j*}\pi_j^*([B])\in A_*(E(X))$.     

Now since $\Phi(\sum_{1\leq j\leq m}b_j\otimes [V_j]) = 0$, we get \\
\noindent
$0~=~[E(V'_k)].\Phi(\sum_{1\leq j\leq m}b_j\otimes [V_j])
~=~\sum_{1\leq j\leq m}\pi^*(b_j)[E(V_j)].[E(V'_k)]
~=~\pi^*(b_k)[E(V(\sigma_k))]$\\
\noindent
$~=~\pi^*(b_k).\iota_{k*}\pi_k^*([B]).$   

Applying $\pi_*$ and using the projection formula we get\\
$0=\pi_*(\pi^*(b_k).\iota_{k*} \pi_k^*([B]))=
b_k . \pi_*\iota_{k*}\pi^*_k([B]) 
=b_k.\pi_{k*}\pi^*_k([B])=b_k.[B]=b_k.$ \\
\noindent
This contradicts our hypothesis that $b_k\neq 0$. It follows that $\Phi$ is 
injective.  

We now prove surjectivity of $\Phi$. One has the filtration $B\cong E(Z_m)
\subset \cdots\subset E(Z_1)=E(X)$.  We claim that  $\Phi$ defines 
surjections 
 $$\Phi_i: A^*(B)\otimes A^*(Z_i)\lr A^*(E(Z_i)) $$ 
for each $i, 1\leq i\leq m$. We prove this by downward induction on $i$. 
This is trivially true for $i=m$, since in this case $E(Z_m)\cong B$.
Consider the diagram \\
$$\begin{array}{cccccc}
A^*(B)\otimes A^*(Z_{i+1})& \lr & A*(B)\otimes A^*(Z_i)&\lr & A^*(B)\otimes 
A^*(Y_i) &\lr 0 \cr
   \Phi_{i+1}~\downarrow~&~&\Phi_i\downarrow &~&~\downarrow ~\cr
A^*(E(Z_{i+1}))&\lr &A^*(E(Z_i)) &  \lr & A^*(E(Y_i)) &\lr 0\\
\end{array} $$
where the top horizontal row is obtained from tensoring with 
$A^*(B)$ the exact sequence $A^*(Z_{i+1}) \lr A^*(Z_i)\lr A^*(Y_i)\lr 0$. 
The homomorphism $A^*(B)\otimes A^*(Y_i)\lr A^*(E(Y_i))$ is 
an isomorphism by prop. 1.9, ch. 1, \cite{f2}. Therefore the 
surjectivity of $\Phi_i$ follows by a diagram chase. \hfill $\Box$

\brem \label{sing}  It follows from the proof of the above lemma that 
the classes  $c_1(\cl_j)\in H^2(E(X);\bz), 1\leq j\leq d,$ generate 
$H^*(E(X);\bz)$ as an $H^*(B;\bz)$-algebra. Similarly, when $p:E\lr B$ is 
algebraic and $B$ a complete nonsingular variety, then 
$[E(V(\rho_j)]\in A^1(E(X)), 1\leq j\leq d,$ generate $A^*(E(X))$ as an  
algebra over $A^*(B)$.
\erem   

We now turn to the ring structure of $H^*(E(X);\bz)$ and $A^*(E(X))$. 

Recall from \S1 that  the line bundle $L_j$ over $X$ 
admits a $T$-equivariant section $s_j:X\lr L_j$ whose zero locus is 
the divisor $V(\rho_j)$. 

Suppose that  
$v_{j_1},\cdots,v_{j_r}$ does not span a cone in $\Delta$. 
Then $s=(s_{j_1},\cdots,s_{j_r})$ is a nowhere vanishing $T$-invariant  
section of $L_{j_1}\oplus\cdots\oplus L_{j_r}$. By taking 
associated  construction, we see that 
the bundle  $\cl_{j_1}\oplus\cdots\oplus\cl_{j_r}$ admits a nowhere vanishing 
section. This implies that 
$$c_1(\cl_{j_1})\cdots c_1(\cl_{j_r})=0 \eqno{(3)} $$ 
in $ H^{2r}(E(X);\bz)$. 

When $p:E\lr B$ is algebraic with $B$ 
nonsingular, we see that  
$$[E(V(\rho_{j_1}))]\cdots[E(V(\rho_{j_r}))]=0 \eqno{(4)}$$ in the 
Chow ring $A^*(E(X))$. 
  
Now, let $u\in M$ be any element. Consider the $T$-equivariant line bundle 
$L_u$ on $X$ corresponding to the principal divisor 
$\sum_{1\leq j\leq d} \langle u,v_j\rangle V(\rho_j)=div(\chi^{-u}).$ 
Clearly $L_u$ is  
isomorphic as a $T$-equivariant bundle to 
$\prod_{1\leq j\leq d} L_j^{\langle u,v_j\rangle}$ 
as both of these bundles correspond to the same piecewise linear 
function $-u:|\Delta|\lr \br$. (See \cite{f}.) 
Hence $E(L_u) \cong  \prod_{1\leq j\leq d} 
\cl_j^{\langle u,v_j\rangle}.$ 
On the other hand the bundle 
$\cl_u:=E(L_u)=E(\chi^{-u}) $ is isomorphic to 
$\pi^*(\xi_1)^{a_1}\cdots \pi^*(\xi_n)^{a_n}$, 
where $a_i=\langle -u, v_i\rangle$. 
This  yields the following relations: 
$$\sum_{1\leq j\leq d} \langle u,v_j\rangle c_1(\cl_j)-
\sum_{1\leq i\leq n}\langle u,v_i\rangle c_1(\pi^*(\xi_i^\vee))=0
\eqno{(5)}  $$ in $H^2(E(X);\bz)$. 
In the case when $p:E\lr B$ is algebraic and $B$ is nonsingular we obtain, in 
the Chow group  $A^1(E(X))$, 
$$\sum_{1\leq j\leq d} \langle u,v_j\rangle [E(V(\rho_j))]
-\sum_{1\leq i\leq n} \langle u, v_i\rangle c_1(\pi^*(\xi_i^\vee))=0.
\eqno{(6)}$$

\noindent 
{\bf Proof of Theorem \ref{main} (i), (iii):}
We first consider part (iii).
In view of equations 
(4) and (6) above we see that we have a well defined homomorphism of algebras:
$\psi:R(A^*(B),\Delta)\lr A^*(E(X))$ defined 
by $\psi(x_j)=[E(V(\rho_j))], 1\leq j\leq d$.
  
Note that, by remark \ref{sing}, $\psi$ is surjective.
We need only prove that $\psi$ is $1-1$. In view of 
theorem \ref{leray}, $A^*(E(X))$ is a  free 
$A^*(B)$-module with basis $[E(V(\tau_i))]$, 
$1\leq i\leq m$.  It follows from lemma \ref{rmod}(iv) that 
$\psi$ is an isomorphism, completing the proof of \ref{main}(iii).

Proof of part (i) is similar. In view of equations (3), (5) above, 
$x_j\mapsto c_1(\cl_j)$ 
defines a homomorphism $R(H^*(B),\Delta)\lr H^*(E(X);\bz)$ 
which is indeed an isomorphism by \ref{rmod}(iv) and \ref{leray}.
\hfill $\Box$

\brem \label{simp}  
If, instead of $\Delta$ being nonsingular, it is only assumed to be 
simplicial, then the toric variety $X$ is only an orbifold. In this case lemma 
\ref{leray} 
holds provided we replace integral homology by rational 
homology and the Chow group by the rational Chow group throughout.
In this case we note that $[V_j].[V_j']= q_j[V(\sigma_j)]$ for a 
rational number $q_j$ and $[V_j][V_k']=0$ for $j>k$. Computing 
the integral cohomology or Chow ring when the fibre $X$ is only simplicial  
seems to be much more difficult.  When $X$ is a weighted projective 
space Al Amrani \cite{ala} has computed the integral cohomology of $E(X)$  
in a more general setting.
\erem   

\section{K-theory}

In this section we prove parts (ii) and (iv) of the main theorem.  

In view of our assumption in \ref{main} (iv) that both the base space 
$B$ and the fibre $X$ are smooth, the Grothendieck ring $\ck^0(E(X))$ of 
algebraic vector 
bundles may be identified, via the duality isomorphism, with  
Grothendieck ring $\ck_0(E(X))$  of coherent sheaves on $E(X)$.  We shall 
denote either of them by $\ck(E(X)$.
Also if a smooth variety $Y$
has an algebraic cell decomposition the forgetful map 
$\ck^0(Y)\lr K(Y)$ is an isomorphism of rings. 
In particular, this holds when  
$X$ is a complete nonsingular toric variety satisfying property $(*)$ 
(see \S 1). 

Although the $K$ ring of 
a complete toric variety has been studied earlier, we could not find 
in the literature its description in terms of generators and relations. 
We obtain such a description in proposition \ref{kx} below.
When $X$ is the projective space, such a description is due to Adams \cite{ad}.
The case when $X$ is a weighted projective space is more recent, due to 
Al Amrani \cite{ala1}. 
We refer the reader to 
\cite{bv},  \cite{m} for other descriptions of the $K$ ring as well as 
the equivariant $K$ ring of a toric variety.  

We begin with the following lemma: 

\blem\label{topline}
Suppose $\zeta_1,\cdots,\zeta_r$ are complex line bundles over 
a finite CW complex $Y$ which has cells only in  even dimensions such that 
$H^*(Y;\bz)$ is generated by 
$c_1(\zeta_1),\cdots, c_1(\zeta_r)\in H^2(Y;\bz)$.  
Then the  ring $K^*(Y)=K^0(Y)$ 
is generated as a ring by $[\zeta_1],\cdots,[\zeta_r]\in K(Y).$   
\elem 
\noindent
{\bf Proof:} Let $f_i:Y\lr \bp^N$ be a classifying map for the bundle 
$\zeta_i, ~1\leq i\leq r$ where $N>1/2(\dim(Y))$. 
Consider the  map $f:Y\lr (\bp^N)^r$ which is defined as 
$f(y)=(f_1(y),\cdots,f_r(y))$.  Then $f^*:H^*((\bp^N)^r;\bz)\lr 
H^*(Y;\bz)$ is easily seen to be a surjection.  
By the naturality of the Atiyah-Hirzebruch \cite{ath} spectral sequence
it follows that $f^*$ induces a surjection of $K$ groups 
$K((\bp^N)^r;\bz)\lr K(Y)$.  Recall from \cite{ad} that 
$K(\bp^N)=\bz[z]/\langle z^{N+1}\rangle $ 
where $z=[\omega]-1$, $\omega$ being the class of the tautological line 
bundle on $\bp^N$.  Hence $K((\bp^N)^r)=\bz[z_1,\cdots,z_r]/\langle 
z_i^{N+1},~1\leq i\leq r\rangle$. 
Since $f^*$ is {\it ring} homomorphism and since 
$f^*(z_i)=[\zeta_i]-1$, the lemma follows. \hfill $\Box$

\blem\label{algline}
Suppose that $Y$ is a complete nonsingular variety over $\bc$ which has 
an algebraic cell decomposition and that $H^*(Y;\bz)$ is generated as 
a ring by $H^2(Y;\bz)$. 
Then there exist algebraic line bundles 
$\zeta_1,\cdots, \zeta_k$ over $Y$ such that $\ck(Y)$  is 
generated as a ring by $[\zeta_i], 1\leq i\leq k$. In particular, 
the forgetful map $\theta:\ck(Y)\lr K(Y)$ is an isomorphism.      
\elem
\noindent
{\bf Proof:} Since $Y$ has an algebraic cell decomposition, 
the Chow ring is 
isomorphic to the cohomology ring $H^*(Y;\bz)$, which is isomorphic as 
an abelian group to $\bz^m$ where $m$ is the number of cells in $Y$.

Since $A^*(Y)$ is torsion free, it follows that $\ck(Y)$ is also torsion free. 
One has the ``topological filtration'' on $\ck(Y)$ and $Gr(\ck(Y))$ denotes 
the associated graded group of $\ck(Y)$. See 15.1.5, \cite{f2}.  
Since the map $\phi: A^*(Y)\lr Gr(\ck(Y))$, defined as $[V]\mapsto [\co_V]$,  
is a surjective homomorphism of groups, it follows that $\ck(Y)$ is a 
finitely generated abelian group of rank at most $m$.  

Let $a_1,\cdots, a_k$ be a $\bz$-basis for $H^2(Y;\bz)$. 
Let $D_1,\cdots, D_k$ be divisors on $Y$ such that $[D_i]\in A^1(Y)$ 
maps to $a_i\in H^2(Y;\bz)$, $1\leq i\leq k$.  
Since the first Chern class of 
$\co(D_i)$ is $a_i\in H^2(Y;\bz)$, for $1\leq i\leq m$, it follows that 
$[\co(D_i)]\in K(Y)$ generate $K^(Y)$ as a ring.  Thus, the forgetful 
homomorphism $\theta: \ck(Y)\lr K(Y)$ is surjective. 

Since the $K(Y)$ is a free abelian group of rank  $m$, 
it follows that $\theta$ is an isomorphism. In particular, 
$\ck(Y)$ is generated as ring by $[\co(D_i)]\in \ck(Y)$, $1\leq i\leq k$.
\hfill $\Box$  

Examples of varieties which satisfy the hypothesis of the 
above lemma are (complete nonsingular) 
toric varieties $X(\Delta)$ where $\Delta$ satisfies $(*)$,  
the flag variety $SL(n,\bc)/B$, where $B\subset SL(n,\bc)$  
is the group of upper triangular matrices, 
and smooth Schubert varieties in $SL(n,\bc)/B$. 

Our next result gives a description of the $K$-ring of $X$. 
We keep the notations of the introduction.

Recall the definition of $\cR$ from \S1. 

\noindent 
\bpropo \label{kx}
Let $X=X(\Delta)$ be a nonsingular complete toric variety where $\Delta$ 
satisfies the property $(*)$. The following relations hold in  $\ck(X)$
and $K(X)$:\\
(i)  $[\co_{V(\rho_{j_1})}]\cdots[\co_{V(\rho_{j_r})}]=0$ if 
$v_{j_1}, \cdots, v_{j_r}$ do not span a cone of $\Delta$,\\
(ii) $\prod_{1\leq j\leq d}[L_j]^{\langle u,v_j\rangle}=1,$\\
(iii) Set $r_i=1\in \bz, ~1\leq i\leq n$.  The homomorphism of rings 
$\theta:\cR=\cR(\bz,\Delta)\lr \ck(X)\cong K(X)$
defnied by $x_i\mapsto [\co_{V(\rho_i)}]=(1-[L_i^\vee])$ is an isomorphism. 
\epropo
\noindent
{\bf Proof:} Recall that $[\co_Y].[\co_Z]=[\co_{Y\cap Z}]$ if $Y,Z$ are 
closed irreducible subvarieties of $X$ which meet transversally. 
Relation (i) follows from the fact that  
$V(\rho_{j_1})\cap\cdots\cap V(\rho_{j_r}) 
=\emptyset $ if $v_{j_1}, \cdots, v_{j_r}$ does not span a cone of $\Delta$. 
Since for any $u\in M$ we have  $\langle -u,v_j\rangle
=\sum_{1\leq p\leq d}\langle u, v_p\rangle \psi_p(v_j)$, 
it follows that one has a $T$-equivariant isomorphism of bundles 
$\prod_{1\leq p\leq d}L_p^{\langle u, v_p\rangle}\cong L_u$, 
where $L_u$ is the line bundle corresponding to the piecewise 
linear function 
$-u:|\Delta|=N_\br \lr \br$. But  $L_u$ is isomorphic to the trivial 
line bundle and so (ii) follows.  

Now the section $s_j:X\lr L_j$ vanishes to order $1$ on $V(\rho_j)$. 
Hence we have an exact sequence of coherent sheaves for $1\leq j\leq d$:
$0\lr L_j^\vee\lr\co_X\lr \co_{V(\rho_j)}\lr 0$. Thus $[\co_{V(\rho_j)}]
=(1-[L_j^\vee])$ in $\ck(X)$, i.e., $(1-[\co_{V(\rho_j)}])=L_j^\vee$. 
Hence $x_j\mapsto [\co_{V(\rho_j}]$ defines a ring homomorphism 
$\theta:\cR\lr \ck(X)$. Since $\ck(X)$ is free abelian 
of rank $m$ and 
since by lemma \ref{kring}  $\cR$ is generated by $m$ elements 
$x(\tau_i), 1\leq i\leq m$,  
it follows that $\theta $ is an isomorphism, completing 
the proof. \hfill $\Box$

\brem\label{sec}
Suppose $\eta_1,\cdots,\eta_k$ are line bundles over $Y$ such that 
their Whitney sum $\eta:=\oplus_{1\leq i\leq k}\eta_i$ admits a nowhere 
vanishing section, then, applying the $\gamma^k$-operation, we obtain 
$\gamma^k(\eta-k)=0$. On the other hand, 
$\gamma^k([\eta]-k)=\gamma^k(\oplus_{1\leq i\leq k}( [\eta_i]-1))
=\prod ([\eta_i]-1)$. Hence, $\prod(1-[\eta_i])=0$. Thus,  
one can avoid the use of the coherent sheaves in the proof of (i) above   
in the case of $K(X(\Delta))$ since 
we know that the section $s=(s_{j_1},\cdots,s_{j_r})$ of the bundle 
$\oplus_{1\leq p\leq r} L_{j_p}$ is nowhere vanishing whenever 
$v_{j_1},\cdots,v_{j_r}$ does not span a cone of $\Delta$.   
\erem

\bcor \label{kxbasis}
(i) The elements $[\co_{V(\tau_i)}]\in \ck(X), 1\leq i\leq m$, form a 
$\bz$-basis for $\ck(X)$. \\
\noindent
(ii) Let  $L(\tau_i)=\prod_{j, v_j\in\tau_i}L_j$ for $1\leq i\leq m$. 
Then $[L(\tau_i)], 1\leq i\leq m,$ form a $\bz$-basis for $K(X)$.
\ecor

\noindent
{\bf Proof:} This follows from the proof of \ref{kx} (iii). \hfill $\Box$ 

Recall that $\cl_j=E(L_j)$ is the line 
bundle over $E(X)$ with total space $E\times_TL_j$. Denote by $\cl(\tau_i)$ 
the line bundle $E(L(\tau_i))=\cl_{j_1}\otimes \cdots \otimes \cl_{j_r}$, 
where $v_{j_1},\cdots, v_{j_r}$ are the primitive vectors along the 
edges of $\tau_i$.
In view of proposition \ref{kxbasis} (ii), the restriction of the 
bundles $\cl(\tau_i),~1\leq i\leq m,$ to the fibre $X$ 
form a $\bz$-basis for $K^*(X)=K^0(X)$.  Hence, by theorem 2.7.8 \cite{at}, 
it follows that $K(E(X))$ is a free $K^*(B)$-module with basis $\cl(\tau_i),
~1\leq i\leq m$. 
Suppose $v_{j_1},\cdots, v_{j_r}$ do not span a cone of $\Delta$. 
The $T$-equivariant section $s=(s_{j_1},\cdots,s_{j_r})$ 
of $L_{j_1}\oplus \cdots \oplus L_{j_r}$ is nowhere vanishing and 
extends to a nowhere vanishing section 
$E(s) :E(X)\lr \cl_{i_1}\oplus\cdots\oplus\cl_{i_r}.$ 
Hence by remark \ref{sec}, 
$$\prod_{1\leq p\leq r}(1-\cl_{j_p})=0.\eqno{(7)}$$ 

Now assume that $p:E\lr B$ is algebraic and $B$ irreducible, nonsingular and 
noetherian over $\bc$. 
Since the $T$-equivariant sections $s_j$ are algebraic, equation (7) 
holds in $\ck(E(X))$ as well.

For any $u\in M$, the $T$-equivariant isomorphism of bundles 
$\prod_{1\leq j\leq d} L_j^{\langle u, v_j\rangle }\cong L_u$ 
yields an isomorphism of vector bundles 
$\prod_{1\leq j\leq d}\cl_j^{\langle u,v_j\rangle}\cong E(L_u).$
Since $E(L_u)=\prod_{1\leq i\leq n}\xi_i^{-\langle u, v_i
\rangle}$, we get 
$$\prod_{1\leq j\leq d} \cl_j^{\langle u,v_j\rangle} 
\cong \pi^*(\xi_{u}^\vee)\eqno{(8)} $$ 
where $\xi_u=\prod_{1\leq i\leq n} \xi_i^{\langle u,v_i\rangle}$.

We are now ready to prove the remaining parts of \ref{main}.

\noindent
{\bf Proof of Theorem \ref{main}(ii), (iv):} 
To prove (ii), we first show that $K^*(E(X))$ is a free $K^*(B)$-module 
of rank $m$, the number of $n$-dimensional cones in $\Delta$ where 
$X=X(\Delta)$. With notations as in \S 4, the restriction of 
$[{\cal L}(\tau_i)], 1\leq i\leq m,$ to the fibre $X$ forms a $\bz$-basis 
for $K^*(X)$. Since $B$ is compact Hausdorff, it is locally 
compact and normal. Therefore $B$ can be covered by finitely many compact 
subsets $W_1,\cdots, W_k$ such that the bundle $\pi|W_r$ 
is trivial for $1\leq r\leq k$. Let $Y$ be any closed subset of $W_r$. 
Now using K\"unneth theorem for K- theory 
(cf. \cite{b}), we see that $K^*(\pi(Y))$ is a free $K^*(Y)$-module 
with basis $[\cl(\tau_i)|Y]$, $1\leq i\leq m$.  Applying Theorem 
1.3, Ch. IV, \cite{k}, we conclude that $K^*(E(X))$ is a free $K(B)$-module 
with basis $[\cl(\tau_i)], ~1\leq i\leq m$.   

In view of equations  (7) and (8), one has a well-defined 
homomorphism of $K(B)$-algebras 
$\Psi:\cR=\cR(K^*(B),\Delta) \lr K^*(E(X))$ 
defined by $x_j\lr (1-[\cl_j]), 1\leq j\leq d.$

Since the $x(\tau_i), 1\leq i\leq m$ span $\cR$ by lemma 
\ref{kring} (iv) and since $K^*(E(X))$ is a free $K^*(B)$ module of 
rank $m$, it follows that $\Psi$ is an isomorphism, completing the 
proof of (ii). 

Now let $B$ be an irreducible, nonsingular, noetherian variety over $\bc$ 
and let 
$p:E\lr B$ be algebraic. Equations (7) and (8) still hold in 
$\ck(E(X))$ since the equivariant 
sections $s_j$ are algebraic.  Proceeding as above, we see that 
to complete the proof of \ref{main} (iv), we need only show that 
$[\co_{E(V_i)}]$, $1\leq i\leq m,$ form a basis for $\ck(E(X))$  
as a $\ck(B)$-module, where $V_i$ stands for $V(\tau_i)$.  
Let $\Phi:\ck(B)\otimes \ck(X)\lr \ck(E(X))$ be the $\ck(B)$-linear map 
defined by $\sum_{1\leq i\leq m}b_i\otimes [\co_{V_i}]\mapsto 
\sum_{1\leq i\leq m}\pi^*(b_i)[\co_{E(V_i)}]$, $1\leq i\leq m$.
In view of \ref{kxbasis}(i), we need only show that 
$\Phi$ is an isomorphism.
 
We first prove surjectivity of $\Phi$. This is proved by induction on the 
dimension of $B$, assuming only that $B$ is noetherian over $\bc$.  
Without loss of generality we may assume that $B$ is irreducible.
If $B$ is a point, then 
the result is obvious. Suppose that $\dim(B)>0$.  Let $U$ be an affine open 
set in $B$ over which the $T$-bundle $p:E\lr B$ is trivial and let 
$Z=B\setminus U$ (with its reduced scheme structure). 
Note that $Z$ may not be irreducible but $\dim (Z_k)<\dim(B)$ for 
each irreducible component $Z_k$ of $Z$. By inductive 
hypothesis, $\ck_0(Z)\otimes \ck(X)\lr \ck(\pi^{-1}(Z))$ is surjective 
homomorphism of abelian groups.  
Consider the commuting  diagram of abelian groups and their homomorphisms:\\
$$\begin{array}{ccccccc}
 \ck_0(Z)\otimes \ck(X)&\lr& \ck_0(B)\otimes\ck(X)&\lr& \ck_0(U)
\otimes\ck(X)&\lr & 0\cr
\downarrow &~& \downarrow &~&\downarrow& ~&~\cr
\ck_0(\pi^{-1}(Z))& \lr&\ck_0(E(X))&\lr &\ck_0(\pi^{-1}(U)) &\lr &0\\
\end{array}$$
where the horizontal rows are exact.  The top horizontal row is got by 
tensoring with $\ck(X)$ the exact sequence 
$\ck_0(Z)\lr \ck_0(B)\lr \ck_0(U)\lr 0$.  
By prop. 2.13, ch. II, (Exp. 0-App.),  p.60, \cite{sga}, the homomorphism 
$\ck(U)\otimes \ck(X)\lr \ck(\pi^{-1}(U))$ is surjective.  
It follows that the homomorphism $\Phi:\ck_0(B)\otimes \ck(X)\lr \ck_0(E(X))$ 
is a surjection.

Now we prove that $\Phi$ is a monomorphism.
Suppose $\Phi(\sum_{1\leq i\leq m}b_i[\co_{V_i}])=0$ where $b_i$ is 
non-zero for some $i$.
Let $p\geq 1$ be the least so that $b_p\neq 0$.  Then, writing 
$V'_p$ for $V(\tau'_p)$, we have  \\
\noindent
$0=[\co_{E(V_p')}].(\sum_{1\leq i\leq m}\pi^*(b_i). 
[\co_{E(V_i)}]) = \sum_{p\leq i\leq m} \pi^*(b_i)[\co_{E(V_p')\cap 
E(V_i))}]$\\ 
\noindent
$=\pi^*(b_p)[\co_{E(V(\sigma_p))}],$ \\
\noindent
since $\Delta$ satisfies property $(*')$. 

Denote by $\pi_p$ the restriction of $\pi$ to $E(V(\sigma_p))$ and 
by $\iota_p$ the  inclusion $E(V(\sigma_p))\subset E(X)$.  Then 
the homomorphism $\iota_{p*}:\ck(E(V(\sigma_p)))\lr\ck(E(X))$ maps 
$[\co_{E(V(\sigma_p))}]=1\in \ck(E(V(\sigma_p)))$ to $[\co_{E(V(\sigma_p))}]
\in \ck(E(X))$. Also, $\pi_p$ is an isomorphism of varieties.
Therefore, applying $\pi_*$ to the expression 
$0=\pi^*(b_p)[\co_{E(V(\sigma_p))}] $ and using the projection formula 
(\S15.1, \cite{f2}) we get \\
$$0=\pi_*(\pi^*(b_p)[\co_{E(V(\sigma_p))}])=b_p.\pi_*\iota_{p*}
([\co_{E(V(\sigma_p))}])=b_p.\pi_{p*}([\co_{E(V(\sigma_p))}])=b_p[\co_B]=b_p.$$
This contradicts our choice of $p$. Hence we conclude that 
$\Phi$ is a monomorphism. \hfill $\Box$

{\bf Concluding remark 4.6.}
Parts (iii) and (iv) of the main theorem also hold when the base field 
$\bc$ is replaced by any algebraically closed field $k$. 
Namely, let $B$ be an irreducible  nonsingular noetherian  
variety over $k$ and 
let $\pi:E\lr B$ be a principal $T$ bundle where $T=Spec(k[M])$.  
Since any toric variety is defined over the integers the fan $\Delta$ defines 
a nonsingular complete $k$-scheme $X=X(\Delta)$.  Again 
$E(X)\lr B$ is a Zariski locally trivial bundle with fibre 
$X$.  Then $\ck(E(X))$ and $A^*(E(X))$ are isomorphic to $\cR$ and 
$R$ respectively.

\end{document}